\documentclass[12pt]{article}
\usepackage{amsmath, amssymb,amsthm}
\usepackage{graphicx}

\newcommand*{\E}{\mathrm{\mathbf{E}}}

\newtheorem{theorem}{Theorem}
\newtheorem{lemma}{Lemma}

\def\eps{\varepsilon}

\def\beqn{\begin{equation}}
 \def\bea{\begin{eqnarray}}
\def\eea{\end{eqnarray}}
\def\bean{\begin{eqnarray*}}
\def\eean{\end{eqnarray*}}

\title{Expansion properties of a random regular
           graph after random vertex deletions}
 \author{
Catherine Greenhill
\thanks{Research supported by the UNSW Faculty Research Grants Scheme.}\\
\small School of Mathematics and Statistics\\[-0.5ex]
\small The University of New South Wales\\[-0.5ex]
\small Sydney NSW 2052, Australia\\[-0.5ex]
\small \tt csg@unsw.edu.au\\
\and
Fred B. Holt\\
\small University of Washington\\[-0.5ex]
\small Seattle WA 98195-4350, USA\\
\small\tt fbholt@u.washington.edu\\
\and
Nicholas Wormald\thanks{Research supported by the Canada Research Chairs program and NSERC.}\\
\small Department of Combinatorics and Optimization\\
\small University of Waterloo\\
\small Waterloo ON, Canada N2L 3G1\\
\small \tt nwormald@uwaterloo.ca\\
}

 \date{}

\begin{document}

\maketitle

\begin{abstract}
We investigate the following vertex percolation process.
Starting with a random regular graph of constant degree,
delete each vertex independently with probability $p$,
where $p=n^{-\alpha}$ and $\alpha=\alpha(n)$ is bounded
away from 0.  We show that a.a.s.\ the resulting
graph has a connected component of size $n-o(n)$ which is
an expander, and all other components are trees of bounded size. 
Sharper results are obtained with extra conditions on $\alpha$. 
These results have an application to the cost of
repairing a certain peer-to-peer network after random failures
of nodes.
\end{abstract}

\section{Introduction}\label{s:intro}

In this paper we investigate the effect of randomly deleting
some vertices in a random regular graph.  
Take a random $d$-regular graph $G$ on $n$ vertices and
independently delete each vertex with probability $p$.
The result is a random graph $\widehat{G}$ with maximum degree
at most  $d$.
We analyse the structure of $\widehat{G}$,
with particular focus on whether (the largest connected
component of) $\widehat{G}$ is an expander graph.  Here $d$ is
fixed, $n$ tends to infinity such that $dn$ is even,
and we take $p=n^{-\alpha}$ for some function $\alpha=\alpha(n)$. 
In this paper we treat only the case where $\alpha$ is bounded away from 0, 
since otherwise even the largest component of the graph is not an expander.
Our work is motivated by an application in peer-to-peer networks,
as described below. 

In Section~\ref{ss:notation} we describe our main result.
Related work is described in Section~\ref{ss:related}.
The application to a certain peer-to-peer network is
explained in Section~\ref{ss:application}.  Our calculations
will be carried out in the configuration model which is
described in Section~\ref{s:configuration}.  Then our calculations
are presented in Section~\ref{s:prob}.  

\subsection{Notation, terminology and our main result}\label{ss:notation}

There are several related definitions of \emph{expander} graphs.
We will say that a
graph $G$ on $n$ vertices is a \emph{$\beta$-expander}
if every set $S$ of $s \leq n/2$ vertices has at
least $\beta s$ neighbours outside $S$.   An alternative definition
involves $d(S)$, the sum of the degrees of vertices in $S$, and
$e(S)$, the number of edges leading out of $S$, and defines $G$ to be 
an $\gamma$-expander if  $e(S)\ge \gamma d(S)$ for all
sets $S\subseteq V(G)$ of vertices with $d(S)\leq |E(G)|$.
For bounded-degree graphs these give equivalent notions of expanders, 
up to a constant factor in translating $\gamma$ to $\beta$.

In this paper, all asymptotics are as $n\to\infty$.  
We say that an event holds
\emph{asymptotically almost surely} (a.a.s.) if the probability
that it holds tends to 1. 
We adapt the standard $O(\cdot)$, $o(\cdot)$ notation to accommodate
versions which hold a.a.s., following~\cite[Section 8.2.1]{NW04}.
Specifically, let $f(n)$, $g(n)$
and $\phi(n)$ be functions such that $|f| < \phi g$.
If $\phi(n)$ is bounded for sufficiently large $n$ then
we write $f=O(g)$, and if $\phi\to 0$ as $n\to\infty$
then we write $f=o(g)$.   When $f/g = 1 + o(1)$ then we
write $f\sim g$ and say that $f$ and $g$ are
\emph{asymptotically equal}.  If a statement $S$
about random variables involves the notations $O(\cdot)$ or
$o(\cdot)$ then $S$ is not an event, and we define ``a.a.s.\ $S$''
to mean that all inequalities of the form $|f| < \phi g$
which are implicit in $S$ hold a.a.s..

Let $\mathcal{G}_{n,d}$ denote the uniform probability space
of all (simple) $d$-regular graphs on the vertex set $[n]=\{ 1,\ldots, n\}$.
Our main result is the following.

\begin{theorem}
\label{main}
Fix $d\geq 3$ and a constant  $\eta >0$.
Suppose that $\alpha=\alpha(n)$ satisfies
\begin{equation}
\label{condition}
\alpha(n) \ge  \eta  
 \end{equation} 
for $n$ sufficiently large.
Let $G\in\mathcal{G}_{n,d}$ and let $\widehat{G}$ be the
graph obtained by independently deleting vertices of $G$
independently with probability $n^{-\alpha}$.  Then \begin{description}
\item{(a)} there is a constant $\beta>0$ such that
 a.a.s.\ $\widehat{G}$ has a connected component of size 
$n-o(n)$ that is a  $\beta$-expander, and all other components are trees
of bounded size;
\item{(b)} 
if $\eta > \frac{1}{2(d-1)}$ then there is a constant $\beta>0$ such that a.a.s.\
$\widehat{G}$ consists of a connected component  that is a  $\beta$-expander, 
together with   $o(n^{(d-2)/(2d-2)})$    isolated vertices;
\item{(c)} if $\eta\ge  {\frac{1}{d-1}}$
then  then there is a constant $\beta>0$ such that a.a.s.\
$\widehat{G}$  is a  $\beta$-expander.
\end{description}
\end{theorem}

The result 
in (a) is best possible, in the sense that if $\alpha$ goes to 0 in the  deletion
probability $n^{-\alpha}$ then there is no fixed positive expansion rate:  that is, 
there is no fixed $\beta >0$ as stated in the theorem. 
The reason for this is as follows. 
It can be shown by the second moment method  that if  
$k<1/(\alpha(d-2))$ then there are a.a.s.\  many paths of 
degree 2 vertices of length at least $k$ in the large connected component. 
Any one such path causes the expansion rate to be at most
at most $2/(k-1)$. This is explained further after Lemma~\ref{step-34} below.

\subsection{Related work}\label{ss:related}

While the vertex-deletion process which we analyse in this paper
does not seem to appear in the
literature, there are various papers~\cite{goerdt, GM, NPSY} 
investigating the result of deleting edges of random regular graphs
independently with some given probability. This is usually described as 
edge percolation, and the resulting graph is
sometimes called the \emph{faulty graph}.
These papers are also motivated by applications to communications networks.
Nikoletseas et al.~\cite{NPSY} focus on the connectivity properties
of the faulty graph, and undertake a study somewhat similar to ours.
Goerdt~\cite[Theorem 2]{goerdt} proves that for small constant edge 
deletion probability,
there is a linear-sized component of the faulty graph. However, it
is \emph{not} an expander.   
Goerdt and Molloy~\cite{GM} extend this analysis to give a threshold
on the fault probability for the existence of a linear sized $k$-core
whenever $3\leq k < d$.  (The $k$-core of a graph is the unique
maximal subgraph in which each vertex has degree at least $k$,
see for example~\cite[p.~150]{bollobas}.) The $k$-core is with high 
probability an expander, but only contains some proportion of the vertices. 
These results are considering much higher deletion probabilities than we do 
in the present paper, because they tolerate a very large number of disconnected 
vertices: linear in $n$.
 
The paper of  Alon et al.~\cite{ABS}  considers edge percolation on expander 
graphs,   which includes random regular graphs of degree at least 3. 
(Though they consider graphs of high girth, this is a minor detail.)
They determine the threshold at which a giant component exists. They also 
give a result~\cite[Proposition 5.1]{ABS} on the expansion of the giant 
component when the edge deletion probability tends to 0. 
This involves $(1/\log n)$ expansion however, not constant rate expansion. 
For random regular graphs, Pittel~\cite{pittel} gave a more detailed analysis 
and determined the order of the transition window of appearance 
of a giant component in a random regular graph under edge percolation.

\subsection{Application to a peer-to-peer network}\label{ss:application}

The vertex deletion process which we study in this paper is
motivated by an application to a 
peer-to-peer network proposed by Bourassa and Holt~\cite{swan-alt-ref,swan}.
This network, called the {\it Swan network},
is based on random regular graphs.  Under normal operating conditions,
the network is given by a $d$-regular graph, where 
$d\ge 4$ is an even constant
(in practice $d=4$).   Bourassa and Holt claimed that their
networks quickly acquire some desirable characteristics of 
uniformly distributed random regular
graphs, such as high connectivity and logarithmic diameter.
(Note that random 
$d$-regular graphs are a.a.s.\ expander graphs for $d\ge 3$~\cite{bollobas88}, 
and as such they are connected and have logarithmic diameter.  Specifically, 
it is 
well known and easy to see that if a graph $G$ is a $\gamma$-expander then 
$G$ has diameter which is bounded above by $\log_{1+\gamma}(n/2)$.)

Cooper, Dyer and Greenhill~\cite{CDG} gave theoretical support
to these claims by defining a Markov chain
to model the behaviour of the Swan networks.   They showed that
under certain natural assumptions about arrival and departure rates, and 
with a slight alteration of the mechanism of departure,
the Markov chain converges rapidly to its stationary distribution, 
which is uniform when conditioned on a fixed number of vertices.  While random 
$d$-regular graphs are a.a.s.\ connected for $d\ge 3$ (indeed, $d$-connected), 
a Swan network in the absence of departures is always connected.

In the context of peer-to-peer Swan networks, the random 
deletion of a vertex corresponds to a client failing.  
Edges correspond to LAN or Internet connections and so are far more 
robust.  Individual clients fail due to lost power, shut down or 
logoff events, frozen applications, and similar phenomena. 
Hence our exclusive consideration of vertex deletions, rather than edge deletions.

Swan networks are self-administering.  In particular, they are self-healing 
after the loss of some vertices, completing a $d$-regular graph among the 
remaining vertices.  For Swan networks, there are two processes to handle 
lost neighbours:  an inexpensive process that uses messages internal to the 
graph, and a more expensive process that contacts vertices using messages 
external to the graph.  As long as the graph remains connected, the repairs 
can safely use the internal repair mechanism. Hence for this application it 
is desirable that the majority of clients in the network remain in a connected 
component. 
 
Theorem~\ref{main} models this situation and shows that the 
large  connected component is an 
expander, which has three important implications for Swan networks.  First, 
under certain constraints on the probability of node failures, Swan networks 
tend to remain connected under the simultaneous loss of several nodes.  
Second, deletions do not degrade the log-diameter of the Swan networks.  
Finally, Theorem~\ref{main} implies that the current internal repair 
strategies could be modified, efficiently involving more of the remaining 
nodes in the repair.

\section{The configuration model and some definitions}\label{s:configuration}

As is usual in this area, calculations are performed in
the \emph{configuration model} (or \emph{pairing model}),
see for example~\cite{NW99} or~\cite[Chapter 9]{JLR}.
A \emph{configuration} consists of $n$ buckets
with $d$ points each, and a perfect matching of the $dn$ points
chosen uniformly at random.  The edges of the perfect matching
are called \emph{pairs}.  Assume that the buckets
are labelled $1,\ldots, n$ and that within each bucket
the points are labelled $1,\ldots, d$.
Denote this probability space by
$\mathcal{P}_{n,d}$.    Given a configuration $P\in\mathcal{P}_{n,d}$
we obtain a pseudograph $G(P)$ by shrinking each bucket down
to a vertex.  This pseudograph may have loops and/or multiple
edges, but the probability that it is simple (with no loops or
multiple edges) is bounded below by a constant.  Moreover,
conditioned on $G(P)$ being simple, it is uniformly distributed.

Similarly if $\mathbf{d} = (d_1,\ldots, d_n)$
is the degree sequence of a graph, then $\mathcal{P}_{n,\mathbf{d}}$
denotes the configuration model where the $j$th bucket contains
$d_j$ points, and a perfect matching of the $2m=\sum_{j=1}^n d_j$
points is chosen uniformly at random.  Here we assume 
that the buckets are labelled $1,\ldots, n$ and that
the points in the $j$th bucket are labelled $1,\ldots, d_j$.

We can now define the bucket deletion process for configurations.
For the remainder of the paper, assume that (\ref{condition}) holds
for some positive constant $\eta$, for $n$ sufficiently large.
Given $P\in\mathcal{P}_{n,d}$, form a new configuration 
$\widehat{P}$ by independently deleting each bucket with probability
$p$.  Specifically:
\begin{itemize}
\item choose a random subset $R$ of buckets such that $b\in R$ with 
probability $p=n^{-\alpha}$, independently for each bucket $b$,
\item delete all buckets in $R$,
\item delete every pair with an endpoint in a bucket in $R$, together with
the other endpoint of the pair if it lies outside $R$,
\item relabel the surviving buckets with the labels $1,2,\ldots$, 
preserving the relative ordering of the buckets,  
\item relabel the points within each surviving bucket in the same way.
\end{itemize}
Note that the same distribution on $\widehat{P}$ will result if
the set $R$ of buckets to delete is chosen first, and then
$P\in\mathcal{P}_{n,d}$ is selected.

We now give some definitions which we will need.
A connected component of a graph which is a tree will
be called an \emph{isolated tree}, and a connected
component of a graph which is a cycle will be called an
\emph{isolated cycle}.

The \emph{$2$-core} of a graph $G$, denoted by $\mathrm{cr}(G)$,
is obtained from $G$
by the following process:  let $G_0=G$ and
for $t\geq 0$, if $G_t$ contains a vertex $v$ of degree 0 or 1
then let $G_{t+1} = G_t - v$, otherwise stop.  The final graph
is $\mathrm{cr}(G)$.  From the 2-core $\mathrm{cr}(G)$ of $G$ 
we obtain the \emph{kernel} of $G$, denoted by
$\mathrm{ker}(G)$, by suppressing all vertices of
degree 2.  That is, if $v$ is a vertex of degree 2 in
$G'$ with neighbours $\{ a, b\}$ then delete $v$ and
replace these two edges by the edge $\{ a,b\}$.

Given a graph $G$, an edge of $G$ is a \emph{cyclic} edge if it belongs
to a cycle, or to a path joining two cycles. 
The cyclic edges are precisely those of the 2-core. The subgraph of $G$
induced by the non-cyclic edges is a union of some number of components.
We call each of these components a \emph{bush}.  If a bush $B$ has
a vertex which is incident with at least one cyclic edge of
$G$ then this vertex is called the \emph{root} of $B$. 
Following from these definitions, a bush can
have at most one root, and the bushes are pairwise disjoint.

We will say that a configuration $P$ has some property if the
corresponding graph $G(P)$ has that property.  This allows
us to speak of paths and cycles in a configuration $P$, 
as well as subconfigurations of $P$ which are trees, bushes
and so on.  In particular we can define the 2-core  and
kernel of a configuration.

We will need the following lemma which 
has a very straightforward proof and can be found in~\cite[p.~54]{bollobas}.
\begin{lemma}
\label{BB-lemma}
Let $k$ be a fixed positive integer and let $\mathbf{d}=(d_1,\ldots, d_n)$
be a degree sequence satisfying $0\leq d_i\leq d$ for all $i$.
Then the probability that a random element of 
$\mathcal{P}_{n,\mathbf{d}}$ contains $k$ specified pairs is
$(1+o(1))(2m)^{-k}$, where $m=(d_1 + \cdots + d_n)/2$ is the number
of pairs in the configuration.
\end{lemma}
 
If an event is a.a.s.\ true for $G(P)$ when $P\in \mathcal{P}_{n,d}$, 
then it is also a.a.s.\ true conditional on the event that $G(P)$ is simple. 
This comes immediately from the fact that the probability that $G(P)$
is simple for $P\in\mathcal{P}_{n,d}$ is bounded below by a nonzero constant 
(see for example~\cite[p.~55]{bollobas}). This is the way that many results   
about $\mathcal{G}_{n,d}$ have been proved using $\mathcal{P}_{n,d}$.

\section{The details}\label{s:prob}

Let $P\in\mathcal{P}_{n,d}$ and let 
$R$ be the random set of buckets chosen
for deletion.   Write $r=|R|$.
By the well known sharp concentration of binomials,   
since $\alpha$ is bounded away from 0,  a.a.s.\ 
\begin{equation}
\label{new}
r\sim n^{1-\alpha} 
\end{equation}
provided $n^{1-\alpha}\to\infty$. Until we come to the proof of 
Theorem~\ref{main} we will assume that the latter condition holds, 
so that (\ref{new}) holds. The other case is easily handled afterwards.

Let $\widehat{P}$ be the result of deleting the
buckets in $R$ from $P$ (and performing the
necessary relabellings of buckets and points).  Then $\widehat{P}$
has $n-r$ buckets.  Let $d_j$ denote the
number of points in the $j$th bucket of $\widehat{P}$, and
say that bucket $j$ has degree $d_j$.
Thus $0\leq d_j\leq d$.  The 
the degree sequence of $\widehat{P}$ is $(d_1,\ldots, d_{n-r})$ and
number of pairs in $\widehat{P}$ is
$(d_1 + \ldots + d_{n-r})/2$. 
Let $N_j$ be the number of buckets of $\widehat{P}$ with degree
$j$, for $0\leq j\leq d$.  
The following result shows that we can use $\mathcal{P}_{n,\mathbf{d}}$ to 
model  $\widehat{P}$, conditional upon it having degree sequence $\mathbf{d}$.

\begin{lemma}
\label{uniform}
The pairing $\widehat{P}$ is uniformly random conditioned on its
degree sequence $\mathbf{d}=(d_1,\ldots, d_{n-r})$.  
\end{lemma}

\begin{proof}
First notice that the set 
$R$ determines an injection $\varphi:[n-r]\to [n]$
which is the inverse of the relabelling operation
performed when $\widehat{P}$ is constructed.
The probability of a particular $\widehat{P}$ with degree sequence 
$\mathbf{d}=(d_1,\ldots, d_{n-r})$ is given by
\[ \binom{n}{{r}} \,\left(\prod_{j=1}^{n-r} \binom{d}{d_j}\right)\, 
                           n^{-\alpha r}\, 
     N_{\widehat{\mathbf{d}}}/|\mathcal{P}_{n,d}|\]
where 
\begin{itemize}
\item $\binom{n}{r}$ is the number of order-preserving injections 
$\varphi:[n-r]\to [n]$, giving the labels of the buckets from $\widehat{P}$
in $P$,
\item $\binom{d}{d_j}$ is the number of order-preserving injections
from $[d_j]$ to $[d]$, giving the labels of the points from
bucket $j$ of $\widehat{P}$ in bucket $\varphi(j)$ of $P$,
\item $n^{-\alpha r}$ is the probability that the $r$ buckets
of $P$ which do not correspond to buckets of $\widehat{P}$ are deleted,
\item $\widehat{\mathbf{d}} = (\widehat{d}_1,\ldots, \widehat{d}_n)$ 
is the degree sequence given by
\[ \widehat{d}_i = \begin{cases} d & \text{if $i\not\in \varphi([n-r])$},\\
                              d-d_j & \text{if $i=\varphi(j)$,}
        \end{cases}\]
\item $N_{\widehat{\mathbf{d}}}$ is the number of configurations
with degree sequence $\widehat{\mathbf{d}}$, 
giving the number of ways to complete the configuration $P$.
\end{itemize}
Since the above expression depends only on $\mathbf{d}$ and not
on the particular structure of $\widehat{P}$, it follows that
$\widehat{P}$ is uniformly random conditioned on its degree sequence
$\mathbf{d}$.
\end{proof}

For $0\leq j\leq d$ let 
\[ \mu_j = \binom{d}{j} n^{1-(d-j)\alpha}.\]

\begin{lemma}
\label{step-1}
Assume that
$r$ satisfies (\ref{new}).  Form $\widehat{P}$ from $P\in\mathcal{P}_{n,d}$
 by deleting the buckets in $R$ as described in Section~\ref{s:configuration}.  
Then,  for $0\leq j\leq d$, we have $\E N_j \sim \mu_j$    and a.a.s.\
\begin{equation}
\label{nj-cases}
 \begin{cases}
    N_j \sim \mu_j  & \text{ if  $\mu_j\to \infty$,} \\
    N_j = O(\log\log n) & \text{ if $\mu_j = O(1)$,}\\
    N_j = 0  & \text{ if $\mu_j = o(1)$.}
    \end{cases}
\end{equation}
In all cases, a.a.s.\ $N_j = o(\mu_{\ell})$ for $0\leq j < \ell\leq d$.
\end{lemma}

\begin{proof}
Fix $j\in \{ 0,\ldots, d\}$.  Choose a random configuration
$P\in\mathcal{P}_{n,d}$. 
The probability that a given bucket $b\not\in R$ is
incident with exactly $d-j$ pairs which are incident
with points in $R$ is asymptotically equal to 
\[ \binom{d}{j}\, \binom{dr}{d-j}\, (d-j)!\, (dn)^{-(d-j)}
      \sim \binom{d}{j} \, n^{-(d-j)\alpha} = \mu_j/n.\]
(The first factor chooses $d-j$ points in $b$ and the
second factor chooses $d-j$ points in $R$.  There
are $(d-j)!$ ways to match up these points using pairs,
and the probability that a random element of $\mathcal{P}_{n,d}$
contains these pairs is $(dn)^{-(d-j)}$, by Lemma~\ref{BB-lemma}.)
Therefore  by linearity of expectation,
\[ \E N_j \sim \mu_j,\]
proving the first statement.

Now suppose that $\mu_j\to\infty$.  Similar calculations for an
ordered pair of buckets $b, c\not\in R$ show that a.a.s.\
$\E [N_j]_2 \sim (\E N_j)^2$.  This establishes the sharp concentration
of $N_j$ whenever $\mu_j\to\infty$.  The other two statements in
(\ref{nj-cases}) follow from Markov's inequality, as does the final
statement of the lemma.
\end{proof}

Now fix a positive integer $K$ such that
\[ K > {\frac{2}{(d-2)\eta}},  \]
where $\eta $ is the constant from (\ref{condition}).
Recall the definition of a bush given before the statement of Lemma 1.
Note that if a bush $B$ in $\widehat{P}$ is not an isolated tree then it has a 
root, and for each non-root bucket $v$ of $B$, the degree of $v$ in $P$ 
is the same as in $\widehat{P}$.

\begin{lemma}
\label{step-2}
Let $\mathbf{d}=(d_1,d_2,\ldots, d_{n-r})$ be a degree sequence
such that $0\leq d_i\leq d$ for all $i$,   $r$ satisfies (\ref{new}), and  
$N_j$ satisfies (\ref{nj-cases}) for $0\le  j\leq d$.
Let $\widehat{P}\in\mathcal{P}_{n,\mathbf{d}}$. 
Then a.a.s.\ $\widehat{P}$ has no bushes
with more than $K$ buckets (including isolated trees). 
\end{lemma}

\begin{proof}
First observe that every tree on $k\geq 2$ buckets has at
least $k/2 + 1$ buckets of degree 1 or 2.
(This can be proved using induction.)
Suppose that $\widehat{P}$ contains a bush $B$
with more than $K+1$ buckets.  Then $\widehat{P}$
contains a bush with exactly $k+1$ buckets, for some $k$ between
$K\leq k\leq 2K$.  To see this, suppose that $B$ has more than
$2K+2$ buckets.  Let $b$ be any bucket of $B$ 
if it is an isolated tree, or let $b$ be the root bucket
of $B$ otherwise.  
Then at least one neighbour of $b$, say $b'$, is the root of a
(smaller) bush $B'$ in $\widehat{P}$ with more than $K$ buckets.  
By induction on $B'$, the result follows.

So now let $B$ be a bush with $k+1$ vertices, where $K\leq k\leq 2K$,
and let $S$ be the set of buckets in $B$.
Ignoring the root bucket (which may have higher
degree in $\widehat{P}$ than it does in $B$),
it follows that there are at least $k/2$ buckets in $S$
with degree 1 or 2 in $\widehat{P}$.  Moreover there are
$k$ pairs in $\widehat{P}$ between points in buckets of $S$.

Now we prove that a.a.s.\ there are no such sets $S$ of
buckets in a random element of $\mathcal{P}_{n,\mathbf{d}}$.
There are $N_1 + N_2$ buckets in $\widehat{P}$ of degree 1
or 2, and by (\ref{nj-cases}), a.a.s.\
\[ N_1 + N_2 = \begin{cases} O(\mu_2) & 
         \text{ if $\mu_2\to\infty$ or $\mu_2=o(1)$,}\\
              O(\log\log n) & \text{ if $\mu_2 = O(1)$.}
              \end{cases}
\]
(Here if $\mu_2\to\infty$ and $\mu_1 = O(1)$ then we use the fact
that a.a.s.\ $N_1 = O(\mu_2)$, rather than the arbitrary upper bound
of $O(\log\log n)$ from (\ref{nj-cases}).)
Hence there are a.a.s.\ at most
\[ O(1)\, n^{k/2+1}\, g(n)^{k/2}\]
ways to choose the buckets belonging to the set $S$, where
\begin{equation}
 g(n) = \begin{cases} {\mu_2} =  n^{(1-(d-2)\alpha)} &
                           \text{ if $\mu_2 \to\infty$ or $\mu_2 =o(1)$,}\\
                \log\log n & \text{ if $\mu_2 = O(1)$.}
          \end{cases}
\label{g-def}
\end{equation}
There are $O(1)$ ways to choose locations for
the $k$ pairs between points of $S$, and the probability
that a random element of $\mathcal{P}_{n,\mathbf{d}}$
contains these pairs is $O(n^{-k})$, by Lemma~\ref{BB-lemma}.  
Therefore the expected number of such sets $S$ in $\widehat{P}$
is a.a.s.\
\[ O(1)\, n^{k/2+1}\, g(n)^{k/2}\, n^{-k}.\]
This is clearly $o(1)$ if $g(n) = \log \log n$, and otherwise 
\begin{align*} 
  O(1)\, n^{k/2+1}\, g(n)^{k/2}\, n^{-k}
  &= O(n^{1-(d-2)\alpha k/2})\\
  &= O(n^{1-(d-2)\eta K/2})\\
  &= o(1)
\end{align*}
by choice of $K$.  Hence by Markov's inequality in either case
there are a.a.s.\ no such sets $S$, for $K\leq k\leq 2K$.
The lemma follows.  
\end{proof}

  To create the 2-core $\mathrm{cr}(\widehat{P})$ of 
$\widehat{P}$,
start with $\widehat{P}$ and delete all buckets of degree 0.
Then while any buckets of degree 1 remain, 
delete one at each time step until none remain.
Finally, relabel the remaining buckets and the points
within the remaining buckets, respecting the relative
ordering.
This process is equivalent to deleting all isolated trees
and ``pruning'' all bushes of $\widehat{P}$ 
(where pruning involves deleting all buckets of the bush 
except the root, and deleting all pairs incident with any non-root
bucket of the bush), followed by relabelling.
Denote the number of buckets in $\mathrm{cr}(\widehat{P})$
by $t$ and let $\mathbf{d}' = (d_1',\ldots, d_t')$ be the
degree sequence of $\mathrm{cr}(\widehat{P})$.  This defines 
$N'_j$, the number
of buckets in $\mathrm{cr}(\widehat{P})$ with degree $j$, 
for $2\leq j\leq d$
(since $\mathrm{cr}(\widehat{P})$ has no buckets of 
degree 0 or 1).

\begin{lemma}
\label{step-34}
Let $\mathbf{d}=(d_1,\ldots, d_r)$ be as in Lemma~\ref{step-2}. 
Let $\widehat{P}\in\mathcal{P}_{n,\mathbf{d}}$.
Then the 2-core $\mathrm{cr}(\widehat{P})$ of $\widehat{P}$ 
has the following properties:
\begin{enumerate}
\item[(i)]  a.a.s.\ $t\sim n-r$ and
$N'_j\sim N_j$ for $2\leq j\leq d$,
\item[(ii)] $\mathrm{cr}(\widehat{P})$ is uniformly random 
   conditioned on its degree sequence,
\item[(iii)] a.a.s.\ $\mathrm{cr}(\widehat{P})$ has no isolated 
cycles,
\item[(iv)] a.a.s.\ $\mathrm{cr}(\widehat{P})$ has no paths
of length at least $K+1$ where all internal vertices have degree 2.
\end{enumerate}
\end{lemma}

\begin{proof}
By Lemma~\ref{step-2}, a.a.s.\ all bushes and
isolated trees in $\widehat{P}$ have at most $K$ buckets
(including the root).  Hence the total number of buckets of
$\widehat{P}$ contained in bushes is a.a.s.\ $O(\mu_1)$ unless
$\lim_{n\to\infty} \mu_1 = O(1)$ and $\mu_1\neq o(1)$, 
in which case an upper bound is given by $O(\log \log n)$.  
Note also that by Lemma~\ref{step-1}, a.a.s.\
\[ N_1 = o(\mu_d) \mbox{ and } \mu_d \sim n-r.\]
It follows that $\mathrm{cr}(\widehat{P})$ has $t$ buckets
where a.a.s.\
\[ t = n-r - o(n-r) \sim n-r.\]
By Lemma~\ref{step-1} again it follows that a.a.s.\
$N'_j\sim N_j$ for $2\leq j\leq d$.
This proves (i).

The proof of (ii) is similar to the argument given in the
proof of Lemma~\ref{step-1} and for similar statements in 
papers on cores of random graphs, so we do not include it here.

Let $m=(d'_1 + \cdots + d'_t)/2$ be the number of pairs
in $\mathrm{cr}(\widehat{P})$.  
By (ii) we know that, conditioned on having degree sequence $\mathbf{d}'$, 
$\mathrm{cr}(\widehat{P})$ has the
distribution of $\mathcal{P}_{n,\mathbf{d}'}$.
Using this and Lemma~\ref{BB-lemma}, 
the expected number of isolated $k$-cycles 
in $\mathrm{cr}(\widehat{P})$ is at most
\[  O(1)\, \binom{N'_2}{k} \, (k-1)!\, (2m)^{-k} =
       O(1)\, \left( {\frac{N'_2}{m}}\right)^k\]
for $2\leq k\leq t$.   
Therefore the expected number
of isolated cycles in $\mathrm{cr}(\widehat{P})$ is at most
\[ O(1) \sum_{k\geq 2} \left({\frac{N'_2}{m}}\right)^k  =
   O(1)\, \left({\frac{N'_2}{m}}\right)^2\, {\frac{1}{1-N'_2/m}} \]
and since a.a.s.\ $N'_2 = o(m)$ by Lemma~\ref{step-1}, we see that
a.a.s.\ the expected number of isolated cycles in $\mathrm{cr}(\widehat{P})$
is $o(1)$.
This establishes (iii), by Markov's inequality.

Finally, the expected number of paths in 
$\mathrm{cr}(\widehat{P})$ of
length $K+1$ with $K$ internal buckets of degree 2 is at most
\[ O(1)\, t^2\, \binom{N'_2}{K}\, K!\, (2m)^{-(K+1)} 
            = O(n)\, \left({\frac{N'_2}{n}}\right)^K\]
using Lemma~\ref{BB-lemma} and (i).
Using (\ref{nj-cases}), a.a.s.\ this expression is
\[ O(n^{1-K})\, g(n)^K \]
where $g(n)$ is defined in (\ref{g-def}).  Using calculations
as in Lemma~\ref{step-2}, this bound is $o(1)$.
Applying Markov's inequality establishes (iv) and completes the proof.
\end{proof}

Following the calculations in this proof, we can now see why
the result of Theorem~\ref{main} is best possible, in the sense
outlined in the introduction.
Suppose that $p=n^{-\varepsilon}$ where $\varepsilon>0$
may be arbitrarily small.  Choose a positive integer $k\geq 2$
such that $\varepsilon (d-2) k < 1$.  (By choosing small enough $\varepsilon$
we may choose $k$ to be arbitrarily large.)
With this deletion probability we have
\[ \mu_2 > \binom{d}{2} n^{1-1/k} \to \infty,\]
so the expected number of paths 
in $\mathrm{cr}(\widehat{P})$ with length at least $k$ and with
at least $k-1$ internal vertices of degree 2 is at least
\[
\Omega(1)\, {\frac{(t-k)^2}{2m}}\, \left({\frac{N'_2 - k}{2m}}\right)^{k-1}
  \geq \Omega(n)\, \left( {\frac{d-1}{2\, n^{1/k}}} \right)^{k-1}\\
  = \Omega(n^{1/k})
  \]
which tends to infinity.   Standard variance calculations show
that the number of such paths is sharply concentrated, so there is
a.a.s.\ at least one such path 
in $\mathrm{cr}(\widehat{P})$.   This implies that the expansion constant
of $\mathrm{cr}(\widehat{P})$  is at most $2/(k-1)$. 
Conditioning on the event that $G(P)$ is simple, we  have the same conclusion 
(see the end of Section~\ref{s:configuration}). 
Hence when $p=n^{-\varepsilon}$ where $\varepsilon = o(1)$, 
we may take $k\to\infty$, there
is no fixed positive expansion rate $\beta$, and the conclusion of 
Theorem~\ref{main} does not hold. 

For practical applications such as the Swan networks, a constant but very small 
deletion probability is the most natural assumption. For the range of 
$n$ of interest in the applications, the probability would be at most  
$n^{-\varepsilon}$ for some small positive constant $\varepsilon$  that is not 
extremely small. For values of the parameters determined in this way, 
we would expect the asymptotic trends studied in this paper to be accurate.

\medskip

\begin{proof}[Proof of Theorem~\ref{main}.] 
Fix a positive integer $d\geq 3$ and a constant  
$\eta >0$
such that (\ref{condition}) holds
for $n$ sufficiently large.
Let $P\in\mathcal{P}_{n,d}$ and form $\widehat{P}$
from $P$ as described in Lemma~\ref{step-1}.  

We first treat the case that $n^{1-\alpha}\to\infty$, so that
(\ref{new}) 
holds a.a.s., and we prove the conclusions of the theorem for the 
multigraph $G(P)$. Only at the end do we remove this assumption and translate 
the result to $G\in \mathcal{G}_{n,d}$.  We have by Lemma~\ref{uniform},
$\widehat{P}\in\mathcal{P}_{n,\mathbf{d}}$ where $\mathbf{d}$
is the degree sequence after deletion.  
Let $\mathrm{cr}(\widehat{P})$ be the 2-core of $\widehat{P}$.
Then a.a.s.\
the conclusion of Lemmas~\ref{step-1},~\ref{step-2},~\ref{step-34} all hold.

Condition on the event that all these conclusions hold, and 
let $\mathrm{ker}(\widehat{P})$ be the kernel of $\widehat{P}$.
Then $\mathrm{ker}(\widehat{P})$ is obtained from $\mathrm{cr}(\widehat{P})$
by suppressing the degree-2 buckets.
That is, if $b=\{ p_1,p_2\}$ is a degree-2 bucket
in $\mathrm{cr}(\widehat{P})$
involved in pairs $\{ p_1,x\}$, $\{ p_2,y\}$, then
delete $b$, remove these pairs and add the pair $\{ x,y\}$.
Since $\mathrm{cr}(\widehat{P})$ has
no isolated cycles, $\mathrm{ker}(\widehat{P})$ 
has exactly $N'_j$ buckets of degree $j$ for $3\leq j\leq d$ 
(and no buckets of degree less than 3).  
For the reasons given in Lemma~\ref{step-34} (ii), we
omit the arguments that show that $\mathrm{ker}(\widehat{P})$ is 
uniformly random conditioned on its degree sequence.
Let $H=G(\mathrm{ker}(\widehat{P}))$ be the multigraph obtained
from $\mathrm{ker}(\widehat{P})$ by shrinking buckets to vertices 
and replacing pairs by edges.   
From~\cite[Lemma 5.3]{BKW}, for some constant $\delta>0$ the multigraph $H$ is 
a.a.s.\ a $\delta$-expander.  (This is well known: for example,  a version of 
this is mentioned in~\cite{goerdt} without proof.)
The constant $\delta$ depends only on $\eta$. 
At this point we further condition on this asymptotically almost sure 
expansion event holding.

Let $G(\widehat{P})$ be the multigraph corresponding to the
pairing $\widehat{P}$.  We obtain $G(\widehat{P})$ from $H$
by performing the following steps:
\begin{itemize}
\item replace some edges by paths of length at most $K$,
\item glue on some bushes of size at most $K$ by identifying their roots with 
       distinct vertices, 
\item introduce some isolated trees of size at most $K$,
\item perform the appropriate relabellings of vertices.  
\end{itemize}
Since we are conditioning on the event that the 
conclusions of Lemma~\ref{step-1} hold,
$G(\widehat{P})$ consists of $O(N_1)=o(n)$ vertices in isolated trees 
of size at most $K$, 
together with a large component having $n-O(N_1)$ vertices. 
Let $U$ be the large component, 
and let $u=|U|$.  Note that $|H|=u-o(u)$.
We now show that $U$ is an expander.

Fix any subset $S\subseteq V(U)$ with $|S|\leq u/2$.  
By an \emph{object} we mean any bush which has been added to $H$, 
or path replacing an edge of $H$, or edge of $H$ not replaced by a path, 
in the process of creating $U$ from $H$. 
An object includes the vertex or vertices of $H$ where it is attached.
An object is \emph{partially occupied} if it has some vertices
in $S$ and some not in $S$, and it is \emph{fully occupied} if all
its vertices belong to $S$.

First suppose that there are at most $\varepsilon|S|/K$ partially
occupied objects, 
where $\varepsilon >0$ is a constant.
Then at most $\varepsilon |S|$ vertices of $S$
are in partially occupied objects.   For each fully occupied object 
there are at most $K-1$ vertices not in $H$ and at least one vertex in $H$. 
Each of these vertices is involved in at most $d$ objects, so the number 
of vertices in $V(H)\cap S$ is at least $1/(d(K-1)+1)> 1/dK$ times the 
number of vertices in fully occupied objects.   
Since all vertices of $S$ are in either partially or fully occupied 
objects, it now follows that
 $$
 |V(H)\cap S| \ge {\frac{1-\eps}{dK}}\, |S|. 
$$

Let $A$ be the set of vertices in $V(H)\setminus S$ which 
have a neighbour in $V(H)\cap S$. We claim that
there exists a constant $\gamma >0$ such that $|A|\ge\gamma\, |V(H)\cap S|$. 
If $|V(H)\cap S| \leq |H|/2$ then the claim follows immediately 
because $H$ is a $\delta$-expander.
So we may assume that $|V(H)\cap S| > |H|/2$. Then $|B|< |H|/2$, 
where 
$B=V(H) \setminus (S\cup A)$, so the expansion of $H$ implies that
$|A|\ge \delta |B|$ and hence 
$$
|A|\ge \frac{\delta}{1+\delta}\, |A\cup B|
=\frac{\delta}{1+\delta} \, |V(H)\setminus  S|
\ge \frac{\delta(1/2+o(1))}{1+\delta}\, |H| 
$$
as $|H|\sim u$ and $|S|\le u/2$. Thus, 
the claim holds with $\gamma= \delta/(3+3\delta)$.
The claim implies that there are at least 
\[ \gamma \, |V(H)\cap S|\ge   {\frac{\gamma(1-\eps)}{dK}}\, |S|\]
partially occupied objects. 

So we may suppose that for some $\eps >0$, there are more than $\eps |S|/K$ 
partially occupied objects. 
Each partially occupied object contains an element of $S$ with a neighbour
in $V(U)\setminus S$.   Since each vertex in $U$ has degree at most
$d$, each of these neighbours can be incident with at most $d$
partially occupied objects.  Therefore $S$ has at least $\eps |S|/dK$
neighbours outside $S$.  It follows from this that the large component $U$
is a constant rate expander, under our assumptions.  The conclusion of part 
(a) of the theorem now follows for the initial random multigraph $G(P)$ in 
place of $G\in\mathcal{G}_{n,d}$, 
under the assumption that $n^{1-\alpha}\to\infty$.

For (b), we need to show further that  a.a.s.\
the only isolated trees in $\widehat{P}$ are isolated vertices.
By Lemma~\ref{step-1}, 
the expected number of isolated trees with
at least two leaves and $k-2$ other vertices is  
\[ n^{k-2}\, n^{2(1-(d-1)\alpha)} \, O\left(n^{-(k-1)}\right)
         = O(n^{1-2(d-1)\alpha}) = o(1).\]
Hence the  isolated trees are a.a.s.\ isolated
vertices, as required.  
Also, the number of isolated vertices is $N_0$
and if $\mu_0\to\infty$ then a.a.s.\
\[ N_0 \sim \mu_0 = n^{1-\alpha d} =o( n^{1-d/(2d-2)})\]
for the given bound on $\eta$.  In the other cases
  we still have $N_0 = o(n^{(d-2)/(2d-2)})$ a.a.s., using Lemma~\ref{step-1}.

For (c), the conclusion of (b) still applies, but in addition, in this 
case $\E N_1 = o(1)$. So there are a.a.s.\ no isolated trees or
bushes of any size, and the conclusion of (c) follows for the multigraph 
$G(\widehat{P})$.  

This completes the proof of the theorem except for two aspects. First, we   
transfer the conclusions from the initial random multigraph $G(P)$ to 
$G\in\mathcal{G}_{n,d}$. This is done  
by conditioning on the event that $G(P)$ is simple.  As explained at the end of 
Section~\ref{s:configuration}, the truth of these asymptotically almost 
sure results is not affected.  
In the conditional space, $G(\widehat{P})$ becomes $\widehat{G}$ as
in the statement of Theorem~\ref{main}.

Finally, we only need to dispense with the assumption that $n^{1-\alpha}\to\infty$. 
Assume  that $n^{1-\alpha}=O(1)$. 
We may apply the version of (c) already proved,   to conclude that 
the deletion of vertices from $G\in\mathcal{G}_{n,d}$ with probability 
$n^{-3/4}$ a.a.s.\ produces a $\beta$-expander $G'$. If we then reinstate   
each deleted vertex (and incident edges) independently with probability 
$1-n^{3/4-\alpha }$ (noting this is positive for $n$ sufficiently large) 
then the result is the same as deleting each vertex of the original graph 
with probability $n^{-\alpha}$, 
that is, it produces $\widehat{G}$.  The vertices deleted from $G$ to 
produce $G'$ are, by easy first moment considerations,  a.a.s.\ of distance 
at least 3 from each other. In this case, reinstating them cannot create any 
new components, and it is easy to see that after 
reinstating them, the resulting graph $\widehat{G}$ is a $(\beta/2)$-expander
when $\beta\leq 2$.   The only nontrivial case is when $S\subseteq V(\widehat{G})$
satisfies $|S - W| < |S\cap W|$, where $W$ is the set of reinstated
vertices.  Here each vertex of $S\cap W$ has $d$ neighbours outside $W$,
giving $d|S\cap W|$ distinct neighbours of $S\cap W$ outside $W$.
At most $|S-W|$ of these can lie in $S$, so $S$ has at least
\[ d|S\cap W| - |S-W| \geq (d-1)|S\cap W|\geq \frac{d-1}{2}|S\cap W|\]
neighbours outside $S$ in $\widehat{G}$.  This gives $\beta/2$-expansion
when $\beta\leq 2$.
\end{proof}

\end{document}